\documentclass[12pt]{article}
\usepackage{amssymb}
\usepackage[top=72pt,bottom=96pt,left=72pt,right=72pt]{geometry}
\def\a{\mathfrak{a}}

\def\e{\varepsilon}
\def\fact[#1,#2]{[\![\,#1\,]\!]_{#2}}
\def\N{\mathbb{N}}
\def\pfill{\par\vskip10pt plus3pt minus3pt\noindent}
\def\proof{\pfill\textbf{Proof:}\quad}
\def\qed{\ensuremath{\hfill\Box}}

\def\R{\mathbb{R}}

\def\Z{\mathbb{Z}}
\newtheorem{Lemma}{Lemma}[section]

\newtheorem{Theorem}[Lemma]{Theorem}

\newtheorem{Corollary}[Lemma]{Corollary}
\newtheorem{Definition}[Lemma]{Definition}

\begin{document}
\title{Sequences of Orthogonal Polynomials related
 to Isotropy Orbits of Symmetric Spaces}
\author{Gregor Weingart\footnote{Unidad Cuernavaca del
 Instituto de Matem\'aticas, Universidad Nacional Aut\'onoma de M\'exico,
 Avenida Universidad s/n, Lomas de Chamilpa, 62210 Cuernavaca, Morelos,
 MEXIQUE;\ \texttt{gw@matcuer.unam.mx}}}
\maketitle
\begin{center}
 \textbf{Abstract}
 \\[15pt]
 \parbox{400pt}{%
  Studying the isotropy orbits of compact symmetric spaces Reiswich
  \cite{reis} introduced a family of explicit polynomials in one
  variable in order to describe the unique minimal isotropy orbit
  of compact symmetric spaces with Dynkin diagram of type $D_m$.
  Based on this geometric interpretation he conjectured that these
  polynomials all have pairwise different real roots in the interval
  $[\,0,1\,]$. In this article the polynomials constructed by Reiswich
  will be identified as special cases of Jacobi polynomials thus proving
  the conjecture about minimal isotropy orbits of compact symmetric
  spaces with Dynkin diagram of type $D_m$.}
 \\[15pt]
\end{center}
\section{Introduction}
 Orthogonal polynomials are a classical topic of study and include many
 famous families of polynomials like the Hermite, Laguerre and Jacobi
 polynomials. Needless to say this article does not strive to make any
 serious contribution to this beautiful topic per se, in fact it will
 turn out eventually that the sequences of orthogonal polynomials studied
 in this article are special cases of Jacobi polynomials. Two of these
 sequences of orthogonal polynomials however appeared recently in a
 parametrization \cite{reis} of the unique minimal isotropy orbit
 of compact symmetric spaces with simply laced Dynkin diagram of type
 $D_m$, more precisely the roots of these polynomials determine the
 coefficients of the minimal isotropy orbit with respect to an orthonormal
 basis of a maximal flat \cite{helg}.

 Mandatory for this geometric interpretation recalled in Corollary \ref{min}
 below is that all roots of the polynomials constructed by Reiswich are
 real, pairwise different and lie in the strict interior of the interval
 $[\,0,1\,]$. All these properties are quite suggestive for a link to
 sequences of orthogonal polynomials, in this article we will establish
 this link and thus prove the conjecture of Reiswich on minimal isotropy
 orbits of a certain family of compact symmetric spaces. The interesting
 question on whether this link extends to general symmetric spaces or not
 remains unanswered for the time being. 
 
 Throughout this article we use the notation $\fact[z,r]$ for the falling
 factorial polynomial instead of the more customary $(\,z\,)_r$, because
 the latter may be too easily confused with other constructs. Recall that
 the falling factorial polynomial is defined for all indices $r\,\in\,\Z$
 by setting $\fact[z,r]\,:=\,z(z-1)\ldots(z-r+1)$ for all positive $r$
 and $\fact[z,0]\,:=\,1$, while $\fact[z,r]\,:=\,0$ in all other cases.
 The falling factorial polynomials make prominent appearance in the
 coefficients of a family of sequences of explicit polynomials
 $(\,R^\tau_n\,)_{n\,\in\,\N_0}$ parametrized by $\tau\,>\,-1$: 

 \begin{Definition}[Family of Reiswich Polynomials]
 \hfill\label{rpol}\break
  For $n\,\in\,\N_0$ the Reiswich polynomial $R^\tau_n\,\in\,\R[\,x\,]$
  with real parameter $\tau\,>\,-1$ is defined by:
  $$
   R^\tau_n(\,x\,)
   \;\;:=\;\;
   \sum_{r\,=\,0}^n(-1)^r\;{n\choose r}\;
   \frac{\fact[n+\tau,r]}{\fact[2n+\tau+1,r]}\;x^{n-r}
  $$
 \end{Definition}

 \noindent
 In his study of orbits of isotropy actions of compact symmetric
 spaces Reiswich \cite{reis} constructed a sequence of polynomials
 $(\,P_m\,)_{m\,\geq\,2}$ in order to characterize the unique minimal
 isotropy orbit for symmetric spaces with simply laced Dynkin diagram
 of type $D_m$. More precisely the original definition of these polynomials
 in \cite{reis} can be rewritten in the form
 \begin{eqnarray*}
  P_m(\,x\,)
  &:=&
  \sum_{r\,=\,0}^{\lfloor\frac m2\rfloor-1}
  (-1)^r\,\left(\;\prod_{d\,=\,1}^r
  \frac{\sum_{\mu\,=\,d+1}^{\lfloor\frac m2\rfloor}(\,1\,+\,2m\,-\,4\mu\,)}
  {\sum_{\mu\,=\,1}^d(\,1\,+\,2m\,-\,4\mu\,)}\;\right)
  \,x^{\lfloor\frac m2\rfloor-1-r}
  \\
  &=&
  \sum_{r\,=\,0}^{\lfloor\frac m2\rfloor-1}
  (-1)^r\;\prod_{d\,=\,1}^r\left(\;\frac{\lfloor\frac m2\rfloor\,-\,d}d
  \;\frac{4m\,-\,4\lfloor\frac m2\rfloor\,-\,4d\,-\,2}
  {4m\,-\,4d\,-\,2}\;\right)\;x^{\lfloor\frac m2\rfloor-1-r}
  \\
  &=&
  \sum_{r\,=\,0}^{\lfloor\frac m2\rfloor-1}
  (-1)^r\,{\lfloor\frac m2\rfloor\,-\,1\choose r}\,
  \frac{\fact[\lceil\frac m2\rceil\,-\,\frac32,r]}{\fact[m\,-\,\frac32,r]}
  \,x^{\lfloor\frac m2\rfloor-1-r}
 \end{eqnarray*}
 where the second line is simply young Gau\ss' formula $a_1+\ldots+a_k
 \,=\,\frac k2(a_1+a_k)$ for arithmetic series. Changing the parameter
 $m\,=\,2n+2$ or $m\,=\,2n+3$ to the parameter $n\,:=\,\lfloor\frac m2
 \rfloor-1$ we conclude $P_{2n+2}(x)\,=\,R^{-\frac12}_n(x)$ and
 $P_{2n+3}(x)\,=\,R^{+\frac12}_n(x)$ respectively.

 \pfill
 In order to study the more general Reiswich polynomials with arbitrary
 real parameter $\tau\,>\,-1$ we consider the probability measure
 $\mu^\tau(dx)\,=\,(\tau+2)\,(\tau+1)\,(1-x)\,x^\tau\,dx$ on the
 interval $[\,0,1\,]\,\subset\,\R$ with an integrable pole at $x\,=\,0$
 for $\tau\,\in\,]\,-1,0\,[$. The moments of the probability measure
 $\mu^\tau(dx)$ are easily calculated via straightforward integration:
 \begin{eqnarray*}
  \mu^\tau_n
  &:=&
  (\tau+2)\,(\tau+1)\,\int_0^1x^n\,(1-x)\,x^\tau\,dx
  \\
  &=&
  (\tau+2)\,(\tau+1)\left.\left(\,\frac{x^{n+\tau+1}}{n+\tau+1}
  \,-\,\frac{x^{n+\tau+2}}{n+\tau+2}\right)\right|_{x=0}^{x=1}
  \;\;=\;\;
  \frac{(\tau+2)\,(\tau+1)}{(n+\tau+2)\,(n+\tau+1)}
 \end{eqnarray*}
 In particular the $0$--th moment $\mu^\tau_0\,=\,1$ tells us that
 $\mu^\tau(dx)$ is in fact a probability measure on $[\,0,1\,]$. In turn
 the probability measure $\mu^\tau(dx)$ gives rise to a positive definite
 scalar product $\langle\,,\rangle$ on the space $\R[\,x\,]$ of polynomials
 with real coefficients by integration against $\mu^\tau(dx)$:

 \begin{Theorem}[Orthogonality Relation]
 \hfill\label{or}\break
  The sequence $(\,R^\tau_n\,)_{n\,\in\,\N_0}$ of Reiswich polynomials is a
  sequence of orthogonal polynomials with respect to the probability measure
  $\mu^\tau(dx)\,:=\,(\tau+2)\,(\tau+1)\,(1-x)\,x^\tau\,dx$ on $[\,0,1\,]$:
  $$
   \int_0^1R^\tau_n(\,x\,)\,R^\tau_m(\,x\,)\,\mu^\tau(\,dx\,)
   \;\;=\;\;
   (n+1)!\,n!\,
   \frac{\fact[n+\tau+1,n]\,\fact[n+\tau,n]}
   {\fact[2n+\tau+2,2n]\,\fact[2n+\tau+1,2n]}\,\delta_{n\,=\,m}
  $$
 \end{Theorem}

 \noindent
 Taking the vanishing orders $\tau$ and $1$ of the weight function $(1-x)
 \,x^\tau$ of the probability measure $\mu^\tau$ at the limits of the
 interval $[\,0,1\,]$ as the decisive clue we obtain as an immediate
 corollary of Theorem \ref{or} that the Reiswich polynomials are special
 cases of Jacobi polynomials $R^\tau_n(x)\,\sim\,P^{(1,\tau)}_n(2x-1)$ up
 to a linear change of variables and normalization \cite{szego}, \cite{wiki}.
 A classical consequence of an orthogonality relation like Theorem \ref{or}
 between the polynomials in a sequence is that the $n$--th Reiswich polynomial
 $R^\tau_n$ with parameter $\tau\,>\,-1$ has exactly $n$ real roots in the
 interior of the support $[\,0,1\,]$ of the probability measure
 $\mu^\tau(dx)$. Combining this statement with the main result of
 Reiswich \cite{reis} we obtain:

 \begin{Corollary}[Minimal Isotropy Orbits \cite{reis}]
 \hfill\label{min}\break
  Recall that the root system of a compact symmetric space $G/K$ equals the
  restriction $\mathfrak{t}^*\longrightarrow\a^*$ of the root system of the
  Lie algebra $\mathfrak{g}$ of the isometry group $G$ with respect to a
  maximal torus $\mathfrak{t}$ to a maximal flat $\a\,\subset\,\mathfrak{t}$.
  For symmetric spaces with Dynkin diagram of type $D_m,\,m\,\geq\,2,$ there
  exists an orthonormal basis $\e_1,\,\ldots,\,\e_m$ of the dual space
  $\a^*$ such that the roots read:
  $$
   D_m
   \;\;\widehat=\;\;
   \{\;\;\pm\,\e_\mu\;\pm\,\e_\nu\;\;|\;\textrm{all choices of signs and\ }
   \;1\,\leq\,\mu\,<\,\nu\,\leq\,m\;\;\}
  $$
  Without loss of generality we may assume that the maximal flat is contained
  $\a\,\subset\,\mathfrak{p}$ in the Cartan complement of the isotropy
  subalgebra $\mathfrak{k}\,\subset\,\mathfrak{g}$ so that the exponential
  map is well--defined by $\exp:\,\a\longrightarrow G/K,\,X\longmapsto
  (\exp\,X)K$. The unique minimal orbit of $K$ on the symmetric space
  $G/K$ passes through the image $\exp(\,X_{\mathrm{min}}\,)K$ under
  $\exp$ of the vector
  $$
   X_{\mathrm{min}}
   \;\;:=\;\;
   \frac\pi4\,E_1
   \;+\;\sum_{r\,=\,1}^n\Big(\;\frac\pi2\,E_{r+1}
   \;+\;\frac{\arccos\sqrt{\xi_r}}2\;(\,E_{m-r}\,-\,E_{r+1}\,)\;\Big)
   \;+\;\frac\pi4\,E_m
  $$
  where $E_1,\,\ldots,\,E_m$ is the orthonormal basis of $\a$ dual to
  the basis $\e_1,\,\ldots,\,\e_m\,\in\,\a^*$ and $0\,<\,\xi_1\,<\,
  \ldots\,<\,\xi_n\,<\,1$ are the $n\,:=\,\lfloor\frac m2\rfloor-1$
  different real zeroes of the polynomial $P_m$.
 \end{Corollary}

 \pfill
 The author would like to express his gratitude for the support and the
 hospitality enjoyed during various stays at the University of Stuttgart,
 in particular he would like to thank U.~Semmelmann and A.~Kollross for
 many interesting mathematical discussions.
\section{The Proof of the Orthogonality Relations}
 In order to prove the orthogonality of the Reiswich polynomials
 $(\,R^\tau_n\,)_{n\,\in\,\N_0}$ formulated in Theorem \ref{or}
 we prove a combinatorial identity, which may be of independent
 interest, between falling factorial polynomials in the first part
 of this section. In a second step we verify that the Reiswich
 polynomials satisfy a recursion formula, whose existence we should
 expect due to orthogonality. Combining both previous results we
 conclude this section with a proof of Theorem \ref{or} and recall
 a standard result on orthogonal polynomials in Corollary \ref{zero}.   

 \begin{Lemma}[Combinatorial Identity]
 \hfill\label{ci}\break
  For every $n\,\in\,\N_0$ and all $u,\,v\,\in\,\N_0$ satisfying $u+v\,
  \leq\,n$ the following polynomial identity holds true in the polynomial
  ring $\Z[\,x,y\,]$ with indeterminates $x,\,y$ and integer coefficients:
  $$
   \sum_{r\,=\,0}^n(-1)^r\,{n\choose r}\,\fact[x-r,u]\,\fact[y-r,v]
   \;\;=\;\;
   n!\,\delta_{u+v\,=\,n}
  $$
 \end{Lemma}

 \proof
 In order to provide the base for an induction on $n\,\in\,\N_0$ let us
 consider $n\,=\,0$ first: The only admissable choice for the additional
 parameters $u,\,v$ equals $u\,=\,0\,=\,v$ in this case leading to the
 trivial identity ${0\choose 0}\,\fact[x,0]\,\fact[y,0]\,=\,1$. Assume now
 by induction hypothesis that the statement of the Lemma holds true for
 $n\,\in\,\N_0$ and all $u,\,v\,\in\,\N_0$ satisfying $u+v\,\leq\,n$.
 Decomposing the binomial ${n+1\choose r}\,=\,{n\choose r}\,+\,{n\choose
 r-1}$ as usual and shifting $r-1$ back to $r$ we obtain for every choice
 $u,\,v\,\in\,\N_0$ of the additional parameters satisfying $u+v\,\leq\,n+1$:
 \begin{eqnarray*}
  \lefteqn{\sum_{r\,=\,0}^{n+1}
  (-1)^r\,{n+1\choose r}\,\fact[x-r,u]\,\fact[y-r,v]}
  \qquad
  &&
  \\
  &=&
  \sum_{r\,=\,0}^n
  (-1)^r\,{n\choose r}\,\Big(\;\fact[x-r,u]\,\fact[y-r,v]
  \;-\;\fact[x-r-1,u]\,\fact[y-r-1,v]\;\Big)
  \\
  &=&
  \sum_{r\,=\,0}^n
  (-1)^r\,{n\choose r}\,\fact[x-r-1,u-1]\,\fact[y-r-1,v-1]
  \,\Big(\,u(y-r)\,+\,v(x-r)\,-\,uv\,\Big)
  \\
  &=&
  +\;u\;\sum_{r\,=\,0}^n(-1)^r\,{n\choose r}\,\fact[x-r-1,u-1]\,\fact[y-r,v]
  \\
  &&
  +\;v\;
  \sum_{r\,=\,0}^n(-1)^r\,{n\choose r}\,\fact[x-r,u]\,\fact[y-r-1,v-1]
  \\
  &&
  -\;u\,v
  \sum_{r\,=\,0}^n(-1)^r\,{n\choose r}\,\fact[x-r-1,u-1]\,\fact[y-r-1,v-1]
 \end{eqnarray*}
 where we have used $(x-r)(y-r)\,-\,(x-r-u)(y-r-v)\,=\,u(y-r)\,+\,v(x-r)
 \,-\,uv$ in the third line; needless to say the cases $u\,=\,0$ and or
 $v\,=\,0$ care for themselves in this calculation. Recall now that
 $u+v\,\leq\,n+1$ holds true by assumption; if this inequality is
 strictly satisfied in the sense $u+v-1\,<\,n$, then the three sums
 on the right hand side all vanish according to our induction hypothesis.
 In the opposite case $u+v\,=\,n+1$ and our induction hypothesis tells
 us that the first two sums both equal $n!$ and the third sum vanishes
 as before so that the right hand side reduces to $(u+v)\,n!\,=\,(n+1)!$.
 \qed

 \begin{Lemma}[Recursion Formulas for Reiswich Polynomials]
 \hfill\label{rec}\break
  Associated to every sequence $(\,R_n\,)_{n\,\in\,\N_0}$ of orthogonal
  polynomials is a recursion formula, which expresses $R_{n+1}$ in terms
  of $R_n$ and $R_{n-1}$ for all $n\,\in\,\N_0$. The recursion formula for
  the sequence of Reiswich polynomials $(\,R^\tau_n\,)_{n\,\in\,\N_0}$ with
  real parameter $\tau\,>\,-1$ reads:
  \begin{eqnarray*}
   R^\tau_{n+1}(\,x\,)
   &=&
   \left(\,x\,-\,\frac{2n^2\,+\,2(\tau+2)n\,+\,(\tau+1)^2}
   {(2n+\tau+3)\,(2n+\tau+1)}\,\right)\,R^\tau_n(\,x\,)
   \\[5pt]
   &&
   \qquad-\;\frac{(n+\tau+1)\,(n+\tau)\,(n+1)\,n}
   {(2n+\tau+2)\,(2n+\tau+1)^2\,(2n+\tau)}\,R^\tau_{n-1}(\,x\,)
  \end{eqnarray*}
 \end{Lemma}

 \proof
 Multiplying the coefficients of the Reiswich polynomials $R^\tau_{n+1}$
 as well as $xR^\tau_n,\,R^\tau_n$ and $R^\tau_{n-1}$ for $x^{n+1-r}$ with
 given index $r\,=\,0,\ldots,n+1$ by the lucky factor
 $$
  C
  \;\;:=\;\;
  (-1)^r\,r\,\frac{2n\,+\,\tau\,+\,1}{{n\choose r-1}}\,
  \frac{\fact[2n+\tau+3,r+1]}{\fact[n+\tau,r-1]}
  \;\;\neq\;\;
  0
 $$
 we obtain by straightforward, if slightly tedious calculation:
 \begin{eqnarray*}
  C\,\Big[\hbox to42pt{\hfill$R^\tau_{n+1}(x)$\hfill}\Big]_{n+1-r}
  &=&
  +\;(n+1)\,(n+\tau+1)\,(2n+\tau+1)\,(2n+\tau-r+3)\,(2n+\tau-r+2)
  \\
  C\,\Big[\hbox to42pt{\hfill$x\,R^\tau_n(x)$\hfill}\Big]_{n+1-r}
  &=&
  +\;(n-r+1)\,(n+\tau-r+1)\,(2n+\tau+3)\,(2n+\tau+2)\,(2n+\tau+1)
  \\
  C\,\Big[\hbox to42pt{\hfill$R^\tau_n(x)$\hfill}\Big]_{n+1-r}
  &=&
  -\;r\,(2n+\tau+3)\,(2n+\tau+2)\,(2n+\tau+1)\,(2n+\tau-r+2)
  \\
  C\,\Big[\hbox to42pt{\hfill$R^\tau_{n-1}(x)$\hfill}\Big]_{n+1-r}
  &=&
  +\;\frac{r\,(r-1)\,(2n+\tau+3)\,(2n+\tau+2)\,(2n+\tau+1)^2\,(2n+\tau)}
  {n\,(n+\tau)}
 \end{eqnarray*}
 Inserting the stipulated coefficients we conclude that the proof of the
 lemma reduces to the finite time fun exercise to verify the following
 degree five polynomial identity in $n,\,\tau$ and $r$:
 \begin{eqnarray*}
  \lefteqn{(n+1)\,(n+\tau+1)\,(2n+\tau+1)\,(2n+\tau-r+3)\,(2n+\tau-r+2)}
  \qquad
  &&
  \\[2pt]
  &=&
  (n-r+1)\,(n+\tau-r+1)\,(2n+\tau+3)\,(2n+\tau+2)\,(2n+\tau+1)
  \\[2pt]
  &&
  \;+\;r\,(2n+\tau+2)\,(2n+\tau-r+2)\,(\,2n^2\,+\,2(\tau+2)n\,+\,(\tau+1)^2\,)
  \\[2pt]
  &&
  \;-\;r\,(r-1)\,(2n+\tau+3)\,(n+\tau+1)\,(n+1)
 \end{eqnarray*}
 \vskip-22pt\qed

 \pfill\textbf{Proof of Theorem \ref{or}:}\quad
 Consider the special case $u\,=\,s$ and $v\,=\,n-s$ of the combinatorial
 identity of Lemma \ref{ci} for some pair of integers $n,\,s\,\in\,\N_0$
 satisfying $s\,\leq\,n$. Evaluating this identity at $x\,=\,n+\tau+s$ and
 $y\,=\,2n+\tau+1$ for a real parameter $\tau\,>\,-1$ leads to
 $$
  n!\;\;=\;\;
  \sum_{r\,=\,0}^n(-1)^r\,{n\choose r}\,
  \fact[n\,+\,\tau\,+\,s\,-\,r,s]\,\fact[2n\,+\,\tau\,-\,r\,+\,1,n-s]
 $$
 Dividing this identity by $\fact[2n+\tau+1,n+1]\,>\,0$ is possible due to
 $\tau\,>\,-1$ and so we find
 \begin{eqnarray*}
  \frac{n!}{\fact[2n\,+\,\tau\,+\,1,n+1]}
  &=&
  \sum_{r\,=\,0}^n(-1)^r\,{n\choose r}\,\frac{\fact[n\,+\,\tau\,+\,s\,-\,r,s]
  \,\fact[2n\,+\,\tau\,-\,r\,+\,1,n-s]}{\fact[2n\,+\,\tau\,+\,1,n+1]}
  \\
  &=&
  \sum_{r\,=\,0}^n(-1)^r\,{n\choose r}\,
  \frac{\fact[n\,+\,\tau,r]}{\fact[2n\,+\,\tau\,+\,1,r]}\,
  \frac1{n\,+\,\tau\,+\,s\,-\,r\,+\,1}
 \end{eqnarray*}
 where $\fact[2n+\tau-r+1,n-s]\,(n+\tau+s-r+1)\,\fact[n+\tau+s-r,s]\,=\,
 \fact[2n+\tau-r+1,n+1]$ should suffice to explain the second equality.
 In the resulting identity the left hand side is independent of $s$,
 subtracting two sucessive instances for $s$, $s+1$ results in the
 key identity
 $$
  \sum_{r\,=\,0}^n(-1)^r\,{n\choose r}\,
  \frac{\fact[n\,+\,\tau,r]}{\fact[2n\,+\,\tau\,+\,1,r]}\;
  \frac1{(\,n+\tau+s-r+2\,)\,(n+\tau+s-r+1\,)}
  \;\;=\;\;
  0
 $$
 valid for all $n,\,s\,\in\,\N_0$ satisfying $s\,<\,n$; the inequality is
 strict now, because we need the previous identity for both $s$ and $s+1$.
 Calculating the scalar product of the Reiswich polynomial $R^\tau_n$ with
 $x^s$ in light of this key identity we obtain directly
 \begin{eqnarray*}
  \langle\,R^\tau_n,\,x^s\,\rangle
  &=&
  \sum_{r\,=\,0}^n
  (-1)^r\,{n\choose r}\,\frac{\fact[n+\tau,r]}{\fact[2n+\tau+1,r]}
  \,\langle\,x^{n-r},\,x^s\,\rangle
  \\
  &=&
  \sum_{r\,=\,0}^n
  (-1)^r\,{n\choose r}\,\frac{\fact[n+\tau,r]}{\fact[2n+\tau+1,r]}
  \;\frac{(\,\tau+2\,)\,(\,\tau+1\,)}{(\,n+\tau+s-r+2\,)\,(n+\tau+s-r+1\,)}
  \;\;=\;\;
  0
 \end{eqnarray*}
 for all $n\,\in\,\N_0$ and $s\,=\,0,\ldots,n-1$ by using the moments
 $\mu^\tau_n\,=\,\frac{(\tau+2)(\tau+1)}{(n+\tau+2)(n+\tau+2)}$ of the measure
 $\mu^\tau(dx)\,=\,(\tau+2)\,(\tau+1)\,(1-x)\,x^\tau\,dx$ calculated before
 in the auxiliary calculation:
 $$
  \langle\,x^{n-r},\,x^s\,\rangle
  \;\;\stackrel!=\;\;
  \mu^\tau_{n+s-r}
  \;\;=\;\;
  \frac{(\,\tau+2\,)\,(\,\tau+1\,)}{(\,n+\tau+s-r+2\,)\,(n+\tau+s-r+1\,)}
 $$
 In consequence $R^\tau_n$ is orthogonal to the subspace of polynomials
 spanned by $1,\,\ldots,\,x^{n-1}$, which evidently contains the Reiswich
 polynomials $R^\tau_0,\,\ldots,\,R^\tau_{n-1}$. Having thus proved the
 orthogonality of the Reiswich polynomials we use the recursion formula
 of Lemma \ref{rec} to calculate their norm squares. Consider the two
 instances of the recursion formula
 \begin{eqnarray}
  R^\tau_{n+1}(\,x\,)
  &=&
  (\,x\,-\,*\,)\;\,R^\tau_n\,(\,x\,)\;\;+\;C\;R^\tau_{n-1}(\,x\,)
  \label{ins1}
  \\
  R^\tau_n(\,x\,)\;\;
  &=&
  (\,x\,-\,*\,)\,R^\tau_{n-1}(\,x\,)\;+\;*\;R^\tau_{n-2}(\,x\,)
  \label{ins2}
 \end{eqnarray}
 for $n$ and $n-1$ where $*$ denotes three irrelevant constants, not
 necessarily the same one. Taking the scalar product of instance
 (\ref{ins1}) with $R^\tau_{n-1}$ we find $0\,=\,\langle\,xR^\tau_n,
 \,R^\tau_{n-1}\,\rangle\,+\,C\,\langle\,R^\tau_{n-1},\,R^\tau_{n-1}
 \,\rangle$, taking similarly the scalar product of instance (\ref{ins2})
 with $R^\tau_n$ we obtain the standard identity
 $$
  \langle\,R^\tau_n,\,R^\tau_n\,\rangle
  \;\;=\;\;
  \langle\,x\,R^\tau_{n-1},\,R^\tau_n\,\rangle
  \;\;=\;\;
  \langle\,R^\tau_{n-1},\,x\,R^\tau_n\,\rangle
  \;\;=\;\;
  -\,C\,\langle\,R^\tau_{n-1},\,R^\tau_{n-1}\,\rangle
 $$
 due to the self--adjointness of the multiplication operator $p\longmapsto
 x\,p$ with the polynomial $x$ under the integration scalar product
 $\langle\,,\rangle$ compare \cite{szego}. Inserting the explicit
 constant $C$ from Lemma \ref{rec} we find eventually the norm square
 of the Reiswich polynomial $R^\tau_n$:
 \begin{eqnarray*}
  \langle\,R^\tau_n,\,R^\tau_n\,\rangle
  &=&
  \frac{(n+\tau+1)\,(n+\tau)\,(n+1)\,n}
  {(2n+\tau+2)\,(2n+\tau+1)^2\,(2n+\tau)}
  \,\langle\,R^\tau_{n-1},\,R^\tau_{n-1}\,\rangle
  \;\;=\;\;\ldots
  \\
  &=&
  (n+1)!\,n!\,\frac{\fact[n+\tau+1,n]\,\fact[n+\tau,n]}
  {\fact[2n+\tau+2,2n]\,\fact[2n+\tau+1,2n]}
  \,\langle\,R^\tau_0,\,R^\tau_0\,\rangle
 \end{eqnarray*}
 \vskip-22pt
 \qed

 \pfill
 For the application of the Reiswich polynomials to the characterization
 of the unique minimal isotropy orbit of a compact symmetric space with
 root diagram of type $D_m$ we recall the following classical statement
 \cite{szego} about polynomials in orthogonal sequences of polynomials:

 \begin{Corollary}[Zeroes of Reiswich Polynomials]
 \hfill\label{zero}\break
  According to a standard result about sequences of orthogonal
  polynomials the $n$--th polynomial $R^\tau_n$ in the Reiswich
  sequence $(\,R^\tau_n\,)_{n\,\in\,\N_0}$ of orthogonal polynomials
  with parameter $\tau\,>\,-1$ has exactly $n$ pairwise different
  real roots in the strict interior $]\,0,1\,[$ of the interval
  $[\,0,1\,]$.
 \end{Corollary}

 \proof
 Polynomial division with remainder provides us in every zero $\xi$
 of a non--zero polynomial $R\,\neq\,0$ with a natural number $o\,\in\,\N$
 and a polynomial $p$ not vanishing in $\xi$ such that $R(\,x\,)\,=\,
 (x-\xi)^o\,p(\,x\,)$. Since $p$ considered as a function is continuous
 and $p(\,\xi\,)\,\neq\,0$, the polynomial $R$ changes its sign in the
 zero $\xi$, exactly if the order $o\,\in\,\N$ of this zero is odd.
 Consider now the set $\{\,\xi_1,\,\ldots,\,\xi_k\,\}$ of all real roots
 of $R^\tau_n\,\neq\,0$ of odd order strictly in the interior of the
 interval $[\,0,1\,]$. In the strict interior of $[\,0,1\,]$ the auxiliary
 polynomial
 $$
  p(\,x\,)\;\;:=\;\;(\,x-\xi_1\,)\,\ldots\,(\,x-\xi_k\,)
 $$
 which reduces to $p(\,x\,)\,=\,1$ in case $k\,=\,0$, changes sign in
 exactly the same points as the polynomial $R^\tau_n$. Replacing $p$ by $-p$
 if necessary we may thus assume that for all $x\,\in\,[\,0,1\,]$:
 $$
  R^\tau_n(\,x\,)\,p(\,x\,)\;\;\geq\;\;0
 $$
 As a non--zero polynomial $R^\tau_n\,p\,\neq\,0$ does not vanish identically
 on any open subset of $[\,0,1\,]$ and so the positivity of the measure
 $\mu^\tau(dx)\,\sim\,(1-x)\,x^\tau\,dx$ ensures the stronger inequality:
 $$
  \langle\,R^\tau_n,\,p\,\rangle
  \;\;:=\;\;
  \int_0^1R^\tau_n(\,x\,)\,p(\,x\,)\,\mu^\tau(\,dx\,)
  \;\;>\;\;
  0
 $$
 On the other hand $R^\tau_n$ is orthogonal to all polynomials of degree
 less than $n$ by assumption
 $$
  R^\tau_n
  \;\;\in\;\;
  \mathrm{span}_\R\{\;R^\tau_0,\,\ldots,\,R^\tau_{n-1}\;\}^\perp
  \;\;=\;\;
  \mathrm{span}_\R\{\;1,\,x,\,\ldots,\,x^{n-1}\;\}^\perp
 $$
 so that the auxiliary polynomial $p$ is necessarily a polynomial of degree
 $k\,\geq\,n$. By construction however $p$ has degree $k\,\leq\,n$ at most
 equal to $n$ so that $R^\tau_n$ has $n$ pairwise different real roots in the
 interior of the interval $[\,0,1\,]$, all of odd degree $o\,=\,1$.
 \qed
\end{document}